\font\ZV=cmr10.pk scaled 2000
.pk scaled 2000
.pk scaled 2000
.pk scaled 2000
.pk
\def\zvezda{\mathop{\,\vrule width0pt depth2pt height8pt
            \smash{\lower7pt\hbox{\ZV *}}\,}\limits}
\font\tenmsam=msam10
\font\sevenmsam=msam7
\font\fivemsam=msam5
\textfont9=\tenmsam
\scriptfont9=\sevenmsam
\scriptscriptfont9=\fivemsam
\mathchardef\le"3936
\mathchardef\ge"393E
\mathchardef\preceq"3934
\mathchardef\succeq"393C

\font\PPPR=msam10
\font\MMMM=msbm10
\def\semitimes{\mathrel{\hbox{\PPPR\char"69}\!}}
\def\N{\hbox{\MMMM N}}
\def\Z{\hbox{\MMMM Z}}


\def\gp#1{\langle #1\rangle}
\def\nc#1{\langle\!\langle#1\rangle\!\rangle}
\let\phi\varphi
\let\~\widetilde

\let \phi \varphi
\let \epsilon\varepsilon
\let \< \langle
\let \> \rangle

\let\=\overline
\def\({{\rm(}}
\def\){{\rm)}}
\def\:{{\rm:}}

\let\s\beginsection
\def\Proof{{\bf\noindent Proof. }}
\frenchspacing
\headline{\ifnum\count0=1 \tt UDC 512.543.7\hss\else\hss\fi}
{\bf
\centerline{HOW TO GENERALIZE THE KNOWN RESULTS}
\centerline{ON EQUATIONS OVER GROUPS}
\smallskip
\centerline{Anton A.~Klyachko}
\it
\centerline{Faculty of mechanics and mathematics,
Moscow State University.}
\centerline{Moscow 119992, Leninskie gory, MSU.}
\centerline{klyachko@daniil.math.msu.su}
}
\footnote{}{This work was supported by the Russian Foundation for 
Basic Research, project no. 05-01-00895.}

\medskip

\s \rm ABSTRACT

The known facts about solvability of 
equations over groups are considered from a more general point of view. A 
generalized version of the theorem about solvability of unimodular 
equations over torsion-free groups is proved. In a special case, this 
generalized version become a multivariable variant of this theorem. For 
unimodular equations over torsion free groups, 
we prove an analogue of Magnus's Freiheitssatz, which asserts that there 
exists a solution with good behavior with respect to free factors of the 
initial group. 

{\it Key words}: equations over groups, Kervaire--Laudenbach conjecture, 
Freiheitssatz. 

\medskip

\s 1. Introduction

An {\it equation over a group $G$ with unknown \({\rm or} variable\) 
$t$} is a formal expression of the form 
$$ 
g_1t^{\epsilon_1}g_2t^{\epsilon_2}\dots g_nt^{\epsilon_n}=1, 
\eqno{(*)} 
$$ 
where $g_i\in G$, $\epsilon_i\in\Z$.  
Equation $(*)$ is called 
{\it solvable over} the group $G$ if there exist a group 
$\~G$ containing $G$ as a subgroup and an element 
$\~t\in\~G$ (called a solution to equation $(*)$) such that
$g_1\~t^{\epsilon_1}g_2\~t^{\epsilon_2}\dots g_n\~t^{\epsilon_n}=1$ in
$\~G$.

There is a lot of theorems (see, e.g., [1]--[14], [18]), wich say that      
equation $(*)$ is solvable under some conditions on the group
$G$ the left-hand side of the equation.
 
In this paper, we suggest a generalization of the notion of equation over 
a group.
{\it A generalized equation over a group $G$ with a variable group $T$} 
is a formal expressio of the form 
$$ 
g_1t_1g_2t_2\dots g_nt_n=1, 
\eqno{(*')} 
$$ 
where $g_i\in G$, $t_i\in T$.  
We call generalized equation $(*')$  
{\it solvable over} the group $G$ if there exist a group 
$\~G$ containing $G$ as a subgroup and a homomorphism 
$T\to\~G,\ t\mapsto\~t$    
(called a solution to generalized equation $(*')$) such that
$g_1\~t_1g_2\~t_2\dots g_n\~t_n=1$ 
in $\~G$.

Clearly, in the case of infinite cyclic variable group, the notion of the 
solvability of a generalized equation coincides with the notion of the 
solvability of an ordinary equation.  
It is suggested to prove generalized analogues of known theorems on 
solvability of equations. Such generalized analogues may include some 
conditions on the variable group, which states that this group is in 
some sense similar to the infinite cyclic.

Let us give some examples. We call an {\it equation}
$$
w(G,t)\equiv \prod g_it^{\epsilon_i}=1,  
\quad\hbox{where }\ g_i\in G,\ \epsilon_i\in\Z,
\eqno{(1)}
$$
or a {\it generalized equation}
$$
w(G,T)\equiv \prod g_it_i=1,
\quad\hbox{where }\ g_i\in G,\ t_i\in T,
\eqno{(1')}
$$
{\it nontrivial} if its left-hand side is not conjugate to a constant 
(i.e., to an element of the initial group $G$) in the free product 
$G*\<t\>_\infty$ (respectively, in $G*T$).

\proclaim{Levin Conjecture} \rm [13].
Any nontrivial equation over a torsion-free group is solvable over this 
group.

This conjecture remains unproven, but it is well-known that it is
equivalent to its generalized analogue.

\proclaim{Generalized Levin conjecture}.
Any nontrivial generalized equation with torsion-free variable group over 
a torsion-free group is solvable over this group.

Indeed, consider a torsion-free group $G_1$ containing the groups $G$
and $T$ as subgroups (e.g., $G_1=G\times T$ or $G_1=G*T$). Starting with 
a nontrivial generalized equation ($1'$) over $G$, we construct the
nontrivial usual equation
$$
v(G_1,t)\equiv w(G_1,t^{-1}Tt)=1
$$
over $G_1$. If the (usual) Levin conjecture is true, then this 
equation has a solution $\~t\in\~{G_1}\supseteq G_1$. But then the group 
$\~T=\~t^{-1}T\~t\subseteq \~{G_1}$ is a solution to the generalized 
equation ($1'$) over the group $G$. This proves the generalized Levin 
conjecture modulo the usual Levin conjecture.

Note that, in the generalized Levin conjecture, the variable group is 
required to be similar to the infinite cyclic in the sense that it 
istorsion-free.

The most known result related to the Levin conjectures 
is the following Brodskii--Howie theorem ([1], [10]). Recall, that 
a group is called {\it locally indicable}, if each of its 
nontrivial finitly generated subgroups admits an epimorphism onto the
infinite cyclic.

\proclaim{Theorem} B. Over a locally indicable group 
any nontrivial equation is solvable.

The generalized analogue of this fact is the following theorem (which
was actually proved in [1]).

\proclaim{Theorem} B$'$. Over a locally indicable group any
nontrivial generalized equation with locally indicable variable
group is solvable.

To prove the equivalence of these two theorems, it is sufficient to 
repeat word-by-word the arguments about the usual and generalized Levin 
conjecture; one should only replace ``torsion-free" by ``locally 
indicable".

Note that, in the Brodskii theorem, the variable group is 
required to be similar to the infinite cyclic in the sense of 
being locally indicable.

We can give many examples of such trivial generalizations. However, it is
not so easy to generalize theorems related to another well-known 
conjecture. Recall, that equation (1) is called {\it nonsingular} if the 
exponent sum of $t$ in $w(G,t)$ is nonzero; if $\sum\epsilon_i=\pm1$, 
then the equation is called {\it unimodular}.

\proclaim Kervaire--Laudenbach conjecture {\rm(see~[3])}. Any unimodular 
equation over any group is solvable over this group.

Sometimes (see, e.g., [5], [10]) the term ``Kervaire--Laudenbach 
conjecture" is used for the hypothesis about solvability of any {\sl 
nonsingular} equation (or even system of equations) over any group. 
At present, no version of this conjecture is proven or disproven.

Trying to formulate a generalized analogue of the 
Kervaire--Laudenbach conjecture or any known related facts, 
we meet the problem: {\sl What should we do with the condition of 
unimodularity (nonsingularity)?}

In this paper we prove a generalized analogue of the following
theorem.

\proclaim {Theorem} 1 \rm{([12], see~also [8])}.
Any unimodular equation over a torsion-free group is solvable over it.

The generalized analogue we are going to prove is as follows.

\proclaim {Theorem} 1$'$.
Any unimodular generalized equation over a torsion-free group is solvable 
over it.

Here, unimodularity of a generalized equation is understood in the sense 
of the following definition that looks not quite natural, but does work 
and, of course, coincides with the usual unimodularity when the 
variable group is infinite cyclic.


\proclaim{Definition} 1.\thinspace\thinspace
\rm Generalized equation ($1'$) is called {\it unimodular} if 
\item{1)} the order of $\prod t_i$ in the group $T$ is infinite; 
\item{2)} $\gp{\prod t_i}\triangleleft T$; 
\item{3)} $T/\gp{\prod t_i}$ is a group with the strong unique 
product property.

Recall that a group $H$ is called a {\it UP-group {\rm or a group with the} 
unique product property} if the product $XY$ of any two finite 
nonempty subsets $X,Y\subseteq H$ contains at least one element, which 
uniquely decomposes into the product of an element from $X$ and an 
element from $Y$. It is known that the group ring of a 
UP-group contains no zero divisors. Some time ago, there was the conjecture 
that any torsion-free group is UP (the converse is, obviously, true).  If 
this conjecture were true, then the Kaplansky problem on zero divisors in 
group rings would be solved. However, it turned out that there exist 
counterexamples ([17], [16]).


We say that a group $H$ has the {\it strong unique product property} if 
the product $XY$ of any two finite 
nonempty subsets $X,Y\subseteq H$ such that $|Y|\ge2$ contains at least 
two uniquely decomposable elements $x_1y_1$ and $x_2y_2$ such that 
$x_1,x_2\in X$,\ \ $y_1,y_2\in Y$ and $y_1\ne y_2$.


Note that the essence of this definition is the condition $y_1\ne y_2$, 
since according to a theorem of Strojnowski [19] the product of any two 
finite nonempty nonsingleton subsets of a UP-group always contains at 
least two uniquely decomposable elements.


As far as we know, all known examples of UP-groups has the strong 
UP-property. In particular, right orderable 
groups, locally indicable groups, diffuse (and weakly diffuse) groups in 
the sense of Bowditch have the strong UP property.


Note also that the strong unique product property follows from the 
following property UP$_4$: the product of any 4 nonempty finite subsets 
$A,B,C,D\subseteq H$ contains at least one element uniquely decomposable 
into a product $abcd$, where $a\in A$, $b\in B$, $c\in C$ and $d\in D$. 
Indeed, if all uniquely decomposable elements $xy$ of the set $XY$ have 
the same factor $y$, then the product of 4 subsets 
$XYY^{-1}X^{-1}$ contains no uniquely decomposable elements: 
the uniqueness decomposability of an element $u=x_1y_1y_2^{-1}x_2^{-1}$ 
implies that $y_1=y$ (otherwise $x_1y_1$ would have another decomposition
$x_1y_1=x'_1y'_1$, and hence $u$ would have another decomposition 
$u=x'_1y'_1y_2^{-1}x_2^{-1}$); the similar arguments show that
$y_2$ shoud be equal to $y$; but then $u$ has two different decompositions
$u=x_1yy^{-1}x_2^{-1}=x_1y'(y')^{-1}x_2^{-1}$, where  $y'\in Y$ is an 
arbitrary element other than~$y$.

\medskip


Theorem $1'$ implies the following multivariable analogue of Theorem 1.

\proclaim{Corollary}.
An equation
$$
g_1x_{j_1}^{\epsilon_1}g_2x_{j_2}^{\epsilon_2}\dots
g_nx_{j_n}^{\epsilon_n}=1
\eqno{(**)}
$$
over a torsion-free group $G$ with variables $x_1,x_2,\dots$ 
is solvable over $G$ if $\prod x_{j_i}^{\epsilon_i}$ is not a proper
power in the free group $F(x_1,x_2,\dots)$.

\noindent{\bf Proof.}
Consider the group
$$
T=\gp{x_1,x_2,\dots\ |\ [x_1,\prod x_{j_i}^{\epsilon_i}]=1,\ 
[x_2,\prod x_{j_i}^{\epsilon_i}]=1,\dots}.
$$

This group is a universal central extension of the one-relator group
$T_1=\gp{x_1,x_2,\dots\ |\ \prod x_{j_i}^{\epsilon_i}=1}$.
It is well-known that, if $\prod x_{j_i}^{\epsilon_i}$ is not a proper 
power in the free group $F(x_1,x_2,\dots)$, then
$T_1$ is locally indicable ([1]) and, hence, $T_1$ is a UP-group. The 
element $\prod x_{j_i}^{\epsilon_i}$ have infinite order in the group 
$T$ 
(see [3]).
Thus, equation $(**)$ can be considered as a 
unimodular generalized equation with variable group $T$. Therefore, 
this equation is solvable by virtue of Theorem  $1'$ and the corollary is 
proven.  

Theorem 1$'$ will be obtained as a special case of the following 
statement.

\proclaim{Main theorem}. If any unimodular equation over
any free power of a group $G$ is Magnusian, then any
unimodular generalized equation over $G$ is solvable over $G$.

We say that equation (1) over a group $G$ is 
{\it magnusian}\footnote{*}
{We mean Magnus's Freiheitssatz [15] (see~also [3]), 
which we can formulate as follows: \sl Any equation over a free group is 
magnusian.} if it is solvable over $G$ and, for any free factor $H$ of 
the group $G$ such that equation (1) is not an equation over it (in the 
sense, that $w$ is not conjugate in $G*\gp{t}$ to an element of the group 
$H*\gp{t}$), some overgroup contains a solution $\~t$ 
transcendental over $H$, i.e., such that 
$\gp{H,\~t}=H*\gp{\~t}_\infty$.

Theorem 1$'$ is a direct corollary of the main theorem and
the following lemma, which is of independent interest.

\proclaim{Lemma} 1.
Any unimodular equation over a torsion-free group is magnusian.

Actually, the main theorem asserts nothing more than the validity of 
Theorem 1$'$, because, obviously, over a group with torsion there exists a 
nonmagnusian unimodular equation (e.g., the equation $t=g$, where $g$ is a 
nontrivial element of finite order).

\proclaim{The notation} \rm which we use is mainly standard. \rm Note 
only that if $x$ and $y$ are elements of a group, and $X$ is a subset 
of this group, then $x^y$ means $y^{-1}xy$, commutator $[x,y]$ is 
understood as $x^{-1}y^{-1}xy$, and the symbols $\gp{X}$ and $\nc{X}$ 
denote, respectively, the subgroup, generated by the set $X$ and the 
normal subgroup, generated by the set $X$. In addition, if $A$ and $B$ are 
isomorphic subgroups of a group and $\phi:A\to B$ is an isomorphism, then 
the system of relations $\{a^x=a^\phi\ |\ a\in A\}$, which often arises in 
our consideration, will be writen in the reduced form $A^x=B$ 
if it is clear what isomorphism is meant.
  

\s 2. Proof of the main theorem

Put $t=\prod t_i$. Let us decompose $T$ into the union of cosets:
$$
T=\coprod_{x\in T/\gp{t}} c_x\gp{t}, \quad\hbox{where }c_1=1.
$$
We rewrite generalized equation ($1'$) in the form
$$
t\prod_i {g_i}^{c_{x_i}t^{k_i}}=1.
\eqno{(2)}
$$
Let $X_1=\{x_i\}$ be the set of all $x\in T/\gp{t}$ occuring in the
irreducible form of equation (2). For each $x\in T/\gp{t}$,
let $G^{(c_x)}$ denote an isomorphic copy of the group $G$ 
assuming that the isomorphism maps $g\in G$ to $g^{(c_x)}\in G^{(c_x)}$.

Put
$$
H_1=\zvezda_{y\in X_1}G^{(c_y)}
$$
and consider the equation 
$$
t\prod_i {g_i}^{(c_{x_i})t^{k_i}}=1
$$
(with variable $t$) over the group $H_1$.

Let us conjugate equation (2) by an element $x\in T/\gp{t}$ (we mean 
that we conjugate by an element of $T$, but conjugation by $t$ gives 
an equation equivalent to the initial one). We obtain the equation 
$$ 
t^{\epsilon_x}\prod_i {g_i}^{{c_{x_ix}}t^{l_i(x)}}=1,
$$
where $\epsilon_x=\pm1$ depending on whether or not $x$ and $t$ commute, 
and the integers $l_i(x)$ are uniquely determined by the equalities 
$c_{x_i}t^{k_i}c_x=c_{x_ix}t^{l_i(x)}$.
Similarly, we put $X_x=X_1x=\{x_ix\}$ and write
the equation 
$$
w_x(t)\equiv t^{\epsilon_x}\prod_i {g_i}^{({c_{x_ix}})t^{l_i(x)}}=1
$$
(with variable $t$) over the group
$$
H_x=\zvezda_{y\in X_x}G^{(c_y)}.
$$
According to Theorem 1, this unimodular equation over a torsion-free group 
has a solution $\~t\in \~{H_x}\supseteq H_x$.

Certainly, we can assume that
$$
\~{H_x}=\<H_x,\ \~t\ |\ w_x(\~t)=1\>.
$$
Put
$$
K=\<\zvezda_{x\in T/\gp{t}}G^{(c_x)},\ \~t\ |\ \{w_y(\~t)=1\ |
\ y\in T/\gp{t}\}\>.
$$
Let us prove that the natural mapping $G\simeq G^{(1)}\to K$ is injective.
To do this, it is sufficient to show that, for any finite set
$Y\subseteq T/\gp{t}$ such that $X_1Y\ni 1$, the natural mapping
$$
G\simeq G^{(1)}\to
K_Y=\<\zvezda_{x\in X_1Y}G^{(c_x)},\ \~t\ |\ \{w_y(\~t)=1\ |
\ y\in Y\}\>
$$
is injective. This is an obvious corollary of the following 
lemma.


\proclaim{Lemma} 2. Suppose that $Y$ is a finite subset of  the group 
$T/\gp{t}$, $z\in T/\gp{t}$,
$k=|X_1|$, and 
$y_1,\dots,y_{k-1}$ are different elements of 
$X_z\cap (X_1Y)$.
Then
in the group $K_Y$ we have
$$
\gp{\~t,G^{(c_{y_1})},\dots,G^{(c_{y_{k-1}})}}=
\gp{\~t}_\infty *G^{(c_{y_1})}*\dots * G^{(c_{y_{k-1}})};
$$
\rm
(Note that the element $z$ may lie either inside or outside $Y$.) 


\Proof
Induction on the cardinality of the set $Y$. If $Y=\{u\}$, then
the assertion of the lemma 
follows from that the unimodular equation $w_u(t)=1$ over the 
group $H_u$ is magnusian (Lemma~1), because $K_Y=K_{\{u\}}=\~{H_u}$.


If $|Y|>1$, then the strong UP property of the group $T/\gp{t}$ implies
that there exists $y\in Y\setminus\{z\}$ such that
$X_y=X_1y\not\subseteq X_1\cdot(\{z\}\cup Y\setminus\{y\})$, therefore
$|(X_1\cdot(\{z\}\cup Y\setminus\{y\}))\cap X_y|<k$; hence, in 
the group $K_{Y\setminus\{y\}}$, by the induction hypothesis we have 
$$
\gp{\~t,\ \{G^{(c_x)}\ ;\
x\in (X_1\cdot(Y\setminus\{y\}))\cap X_y\}}=
\gp{\~t}_\infty*\left(\zvezda_
{(X_1\cdot(Y\setminus\{y\}))\cap X_y}
G^{(c_x)}\right).
$$
Therefore, in the group 
$$
N=K_{Y\setminus\{y\}}*
\left(\zvezda_
{x\in (X_z\cap X_y)\setminus(X_1(Y\setminus\{y\}))}
G^{(c_x)}\right)
$$
we have
$$
M=
\gp{\~t,\ \{G^{(c_x)}\ ;\
x\in (X_1\cdot(\{z\}\cup Y\setminus\{y\}))\cap X_y\}}=
\gp{\~t}_\infty*\left(\zvezda_
{x\in (X_1\cdot(\{z\}\cup Y\setminus\{y\}))\cap X_y}
G^{(c_x)}\right).
$$
The same decomposition of the group $M$ is valid in $\~H_y$ by Lemma~1. 
Thus, we obtain a correct decomposition into the amalgamated product  
$$ 
K_Y=N\zvezda_M\~H_y.  
$$
The subgroup 
$\gp{\~t,G^{(c_{y_1})},\dots,G^{(c_{y_{k-1}})}}$
lies in $N$ and the assertion of the lemma follows from the induction 
hypothesis. 
Lemma 2, and the injectivity of the natural mapping
$G\simeq G^{(1)}\to K$ is proven.

The group $T$ acts on $K$ by the rule
$$
\~t^y=\~t^{\epsilon_y},\quad
(g^{(c_x)})^y=(g^{(c_f)})^{\~t^k}, \eqno{(3)}
$$
where $y\in T$ is an arbitrary element, and $f\in T/\gp{t}$ and $k\in\Z$
are determined from the equalities $c_xy=c_ft^k$. Obviously, this is a 
well defined action by automorphisms.

Take the corresponding semidirect product
$T\semitimes K$ and consider its quotient group by the 
normal cyclic
subgroup
$\gp{\~tt^{-1}}$. This is the sought group containing a solution:
$$
\~G=(T\semitimes K)/\gp{\~tt^{-1}}.
$$
Indeed, $G$ is embedded into $\~G$ as subgroup:
$G=G^{(1)}\subseteq K \subseteq \~G$. According to the definition of action
(3), we have
$$
G^{(c_x)}=G^{c_x},
$$
hence, the relation $w_1(\~t)=1$ holds in $K$ and equality $t=\~t$ in 
$\~G$ give relation (2). Thus, the group $T$  considered as a 
subgroup of $\~G$ is a solution to equation (2). The main theorem 
is proven.


\s 3. Proof of Lemma 1

First, note one simple fact.

\proclaim{Lemma} 3. Suppose that $A$ is a nontrivial subgroup of a group 
$B$ and $b\in B$. Then $b$ is transcendental over $A$ if and only if 
$$ 
\gp{\{A^{b^i}\ | \ i\in \Z\}}=\zvezda_{i\in \Z} A^{b^i}.
$$

\Proof
The ``only if" part is obvious. The ``if" part follows from 
the fact that, if $u\in A*\gp{b}_\infty$ is a nontrivial relation between 
$A$ and $b$ in the group $B$ and $a\in A\setminus\{1\}$, then $[a,u]$ is a 
nontrivial relation between the groups $A^{b^i}$. Lemma 3 is proven.

Throughout this section, we suppose that $G=H*K$, equation (1) 
is an unimodular equation over $G$, and $w$ does not belong to a 
subgroup conjugate to $H*\gp{t}$ in the group $G*\gp{t}_\infty$. In 
addition, we put $U=G*\gp{t}_\infty/\nc{w}$ (this is the so caled 
universal solution group). We have to prove that equation (1) has a 
solution transcendental over $H$, or, equivalently, the element $t\in U$ 
is transcendental over $H$.  

\proclaim{Lemma} 4. If 
$$ 
\gp{\{H^{t^i}\ | \ i\in \Z\}}=\zvezda_{i\in \Z} H^{t^i}
\eqno{(4)} 
$$ 
in $U$, then $t\in U$ is a solution to equation $(1)$ transcendental over 
$H$.

\Proof
In the case when the group $H$ is nontrivial, the assertion immediately 
follows from Lemma 3. If $H=\{1\}$, then the lemma asserts only that 
the element $t\in U$ has infinite order (if $w$ is not conjugate to 
$t^{\pm1}$); this fact is observed in [6]. Lemma 4 is proven.


Put
$$
\=H=(\zvezda_{i\in \Z} H_i),\quad \=G=\=H*K,
$$
where $H_i$ are isomorphic copies of the group $H$. Identifying the group 
$G$ with the subgroup $N_0*K$ of $\=G$ consider the system of equations 
over $\=G$ consisting of initial equation (1) and auxiliary 
equations 
$$ 
\Biggl\{\eqalign{ w(t)&=1,\cr t^{-1}H_it&=H_{i+1},\quad 
i\in\Z.\cr } 
$$ 
Using the obvious transformations this system can be rewriten in the form 
$$ 
\Biggl\{\eqalign{ w_1(t)&=1,\cr t^{-1}H_it&=H_{i+1},\quad i\in\Z,\cr } 
$$ 
where $w_1$ contains no fragments of the form $t^{-1}\=ht$ and $t\=ht^{-1}$ 
with $\=h\in\=H$.


If the length of $w_1(t)$ is 1, then the first equation of this system
can be rewriten in the form $t=u$, where $u\in\=G$, 
and the entire system can be rewriten in the form
$$
\{
u^{-1}H_iu=H_{i+1},\quad i\in\Z\}.
$$

The natural mapping of the group $\=N$ into the group 
$\={\=G}=\gp{\=G\ |\ \{u^{-1}H_iu=H_{i+1},\quad i\in\Z\}}$
is injective by the following theorem proven in [12] (see
also [8]).

\proclaim{Theorem}. Suppose that $A$ and $B$ are torsion-free groups, 
$v\in (A*B)\setminus A$, and $\phi$ is an automorphism of $A$. Then  
the natural mappings
$$
A\to\gp{A*B\ |\ \{a^v=a^\phi\ |\ a\in A\}} \leftarrow B
$$
are injective.

The element $t$ of the HNN-extension 
$\~G=\gp{\={\=G},t\ | H_i^{t}=H_{i+1}\ i\in\Z}$ is obviously a solution 
to equation (1) over $G=H_0*K$. Clearly, it is transcendental over 
$H=H_0$ by Lemma 4, because relation (4) is fulfilled in $\~G$ and, hence, 
in $U$. (Actually, it is easy to observe that $\~G\simeq U$.)

Now, consider the case when the length (that is, the number of occurrences 
of $t^{\pm1}$) of the word $w_1$ is greater than one.
Consider the following subgroups of $G*\gp{t}_\infty$:
$K_i=t^{-i}Kt^i$,
$H_i=t^{-i}Ht^i$,
$$
\=H=\zvezda_{i=-\infty}^\infty H_i,\quad
K^{(m)}=\zvezda_{i=0}^m K_i \quad\hbox{and}\quad G^{(m)}=\=H*K^{(m)}.
\eqno{(5)}
$$
Consider all possible expressions of 
equation (1) in the form
$$
ct\prod_{i=1}^n b_it^{-1}a_it=1, \quad\hbox{where } a_i,b_i,c\in G^{(m)}.
\eqno{(6)}
$$
Among all such expressions we choose those in which $m$ is minimal; after 
that, from all expressions with minimal $m$ we choose an expression with 
minimal $n$.  For such a minimal expression (6), we have:  
\item{1)} 
$n\ge 1$ 
(i.e., the length of this expression is strictly greater than one); 
\item{2)} 
$a_i\notin G^{(m-1)}=\=H*K_0*\dots*K_{m-1}$ and 
$b_i\notin (G^{(m-1)})^t=\=H*K_1*\dots*K_m$; 
\item{3)} 
$a_i$ are transcendental over 
$G^{(m-1)}$, and $b_i$ are transcendental over $(G^{(m-1)})^t$.

The first property is ensured by the assumption that the length of $w_1$ is
greater than one (in any expression of length 1 we have $m>0$; 
therefore, we can decrease $m$ by replacing all occurrences of elements of
$K_m$ by fragments $t^{-1}gt$, where $g\in K_{m-1}$). The second property 
obviously follows from the minimality of $n$ and $m$. Property 3) 
is a corollary of property 2).

Now, suppose that the symbols $H_i$ and $K_i$ denote abstract isomorphic 
copies of the groups $H$ and $K$, and the groups $\=H$, $K^{(m)}$ and 
$G^{(m)}$ are defined by formulae (5). Consider the following system 
of equations over the group $G^{(m)}$:  
$$ 
\left\{\eqalign{ 
&x^{-1}H_ix=H_{i+1},\quad i\in\Z, \cr &x^{-1}K_ix=K_{i+1},\quad 
i\in\{0,\dots m-1\}, \cr &cx\prod_{i=1}^n b_ix^{-1}a_ix=1.  }\right.  
\eqno{(7)}
$$
Clearly, any solution to this system over $G^{(m)}$ is a solution 
to equation (1) over $G$ transcendental over $H$ (by Lemma 4). To 
complete the proof of Lemma 1, we should only note that 
properties 1) and 3) of system (7) imply its solvability in virtue of 
the following theorem. 

\proclaim{Theorem} \rm{([12], see also [8])}. Suppose that $M$ and $N$ are
isomorphic subgroups of a group $L$, $\phi:M\to N$ is an isomorphism, 
$n\in\N$, $a_1,\dots, a_n$ are elements of $L$ transcendental over 
$M$, $b_1,\dots, b_n$ are elements of $L$ transcendental over $N$, and 
$c\in L$. 
Then the system of equations 
$$ 
\left\{\eqalign{ &x^{-1}gx=g^\phi,\quad 
g\in M, \cr &cx\prod_{i=1}^n b_ix^{-1}a_ix=1 }\right.  
$$ 
is solvable over $L$.


\s 4. Concluding remarks

The definition of the unimodularity of generalized equations (Definition 
1) would look nicer, if condition 3) had the form 
\item{$\~3)$} the group $T/\gp{\prod t_i}$ is torsion-free.

\noindent We call generalized equation $(1')$ satisfying conditions
1), 2) and $\~3)$ {\it weakly unimodular}. 

\proclaim Question 1.
Is it true that any weakly unimodular generalized equation
over a torsion-free group is solvable over it?

If we completly remove condition 3) from Definition 1, we shall obtain the 
notion of a {\it nonsingular} generalized equation.

\proclaim Question 2.
Is it true that any nonsingular generalized equation over a torsion-free 
group is solvable over it?

The answer to this question is unknown even for usual equations, in the 
generalized situation it might be easier to construct a counterexample.

Is it possible to prove generalized analogues of others known facts? for 
example, the Gerstenhaber--Rothaus theorem ([9], see~also [3])?

\proclaim Question 3.
Is it true that any \(weakly\) unimodular \(or even any
nonsingular\) generalized equation over a finite group is solvable over 
it?

It is easy to see that the Levin conjecture implies that any solvable
equation over a torsion-free group is magnusian. 
Let us state the following conditional question.

\proclaim Question 4.
Is the conjecture that any solvable 
equation over a torsion-free group is magnusian equivalent to 
the Levin conjecture?

It is easy to see that the answer to the generalized version of this 
question is positive.

\s{\rm REFERENCES}

\frenchspacing
\item{[1]} 
Brodskii S. D.
Equations over groups and one relator groups (in Russian) 
{// Sib. Mat. Zh.} 1984. {T.25}. \number 2. S.84--103.

\item{[2]} 
Klyachko Ant. A., Prishchepov M. I.
The descent method for equations over groups (in Russian)
{// Vestn. MSU: Mat., Mech.} 1995. {\number4}. S.90--93.

\item{[3]}
R. C. Lyndon and P. E. Schupp, 
{Combinatorial group theory}, 
Springer-Verlag, 1977.

\item{[4]} 
Clifford A.
Nonamenable type K equations over groups
{// Glasgow Math. J.} 2003. {V.45}. P.389--400.

\item{[5]} 
Clifford A., Goldstein R. Z.
Equations with torsion-free coefficients
{// Proc. Edinburgh Math. Soc.} 2000. {V.43}. P.295--307.

\item{[6]} 
Cohen M. M., Rourke. C.
The surjectivity problem for one-generator, one-relator extensions of 
torsion-free groups
{// Geometry \& Topology}. 2001. {V.5}. P.127--142.

\item{[7]} 
Edjvet M., Howie J.
The solution of length four equations over groups
{// Trans. Amer. Math. Soc.} 1991. {V.326}. P.345--369.

\item{[8]} 
Fenn R., Rourke C.  
Klyachko's methods and the solution of equations over torsion-free groups
{// L'Enseignment Math\'ematique.} 1996. {T.42}. P.49--74.

\item{[9]} 
Gerstenhaber M., Rothaus O. S.
The solution of sets of equations in groups
{//  Proc. Nat. Acad. Sci. USA}. 1962. {V.48} P.1531--1533.

\item{[10]} 
Howie J.
On pairs of 2-complexes and systems of equations over groups
{// J. Reine Angew Math.} 1981. {V.324}. P.165--174.

\item{[11]} 
Ivanov S. V., Klyachko  Ant. A.
Solving equations of length at most six over torsion-free groups
{// J. Group Theory}. 2000. {V.3}. P.329--337.

\item{[12]} 
Klyachko  Ant. A.
A funny property of a sphere and equations over groups
{// Comm. Algebra}. 1993. {V.21}. P.2555--2575.

\item{[13]} 
Levin F.
Solutions of equations over groups
{// Bull. Amer. Math. Soc.} 1962. {V.68}. P.603--604.

\item{[14]} 
Lyndon R. C.
Equations in groups
{// Bol. Soc. Bras. Math.} 1980. {V.11}. \number1. P.79--102.

\item{[15]} 
Magnus W.
\"Uber diskontinuierliche Gruppen mit einer definierenden Relation
(Der Freiheitssatz).
{// J. Reine Angew Math.} 1930. {V.163}. P.141--165.

\item{[16]} 
Promyslow S. D.
A simple example of a torsion free nonunique product group
{// Bull. London Math. Soc.} 1988. {V.20}. P.302--304.

\item{[17]} 
Rips E., Segev Y.
Torsion free groups without unique product property
{// J. Algebra} 1987. {V.108}. P.116--126.

\item{[18]} 
Stallings J. R. 
A graph-theoretic lemma and group embeddings 
// Combinatorial group theory and topology 
(ed. Gersten S. M, Stallings J. R.). 
{Annals of Mathematical Studies}. 1987. {V.111}. P.145--155.

\item{[19]}
Strojnowski A.
A note on u.p. groups
{// Comm. Algebra} 1980. {V.8}. P.231--234.

\end